\documentclass{amsart}

\theoremstyle{definition}

\theoremstyle{remark}

\begin{document}

\title[A quartic integral]{A formula for a quartic integral: a survey of old
proofs and some new ones}

\author{Tewodros Amdeberhan}
\address{Department of Mathematics,
Tulane University, New Orleans, LA 70118}
\email{tamdeberhan@math.tulane.edu}

\author{Victor H. Moll}
\address{Department of Mathematics,
Tulane University, New Orleans, LA 70118}
\email{vhm@math.tulane.edu}

\subjclass{Primary 33}

\date{\today}


\begin{abstract}
We discuss several existing proofs of 
the value of a quartic integral and present 
a new proof that evolved from rational Landen transformations. We include 
our personal renditions as related to the history of these particular 
computations.
\end{abstract}

\maketitle

\newcommand{\nn}{\nonumber}
\newcommand{\ba}{\begin{eqnarray}}
\newcommand{\ea}{\end{eqnarray}}
\newcommand{\ift}{\int_{0}^{\infty}}
\newcommand{\ifft}{\int_{- \infty}^{\infty}}
\newcommand{\no}{\noindent}
\newcommand{\X}{{\mathbb{X}}}
\newcommand{\Q}{{\mathbb{Q}}}
\newcommand{\R}{{\mathbb{R}}}
\newcommand{\Y}{{\mathbb{Y}}}
\newcommand{\Ftwo}{{{_{2}F_{1}}}}
\newcommand{\realpart}{\mathop{\rm Re}\nolimits}
\newcommand{\imagpart}{\mathop{\rm Im}\nolimits}

\newtheorem{Definition}{\bf Definition}[section]
\newtheorem{Thm}[Definition]{\bf Theorem}
\newtheorem{Example}[Definition]{\bf Example}
\newtheorem{Lem}[Definition]{\bf Lemma}
\newtheorem{Note}[Definition]{\bf Note}
\newtheorem{Cor}[Definition]{\bf Corollary}
\newtheorem{Prop}[Definition]{\bf Proposition}
\newtheorem{Problem}[Definition]{\bf Problem}
\numberwithin{equation}{section}

\section{Introduction} \label{sec-intro}
\setcounter{equation}{0}

The evaluation of definite integrals has attracted the scientific community, 
both professional and amateurs, for a long time.  Some of these evaluations
are collected in tables. The authors' favorite one is by I. Gradshteyn and 
I. Ryzhik \cite{gr}. It was the  naive attempt to prove all the formulas in
\cite{gr} that produced \cite{irrbook}. A recent review of this 
book \cite{foncannon} 
describes how R. Feynman 'once claimed that he acquired his initial 
reputation  not as a physicist, but as a redoubtable evaluator of integrals'. 

In the Fall of 1995, when the author (VHM) was teaching a graduate class 
in Analysis, a student (George Boros) stated that he could prove the identity
\begin{equation}
\ift \frac{dx}{(x^{4} + 2ax^{2} +1 )^{m+1}} = \frac{\pi}{2} 
\frac{P_{m}(a)}{[ 2(a+1) ]^{m + \tfrac{1}{2}}},  
\label{int-01}
\end{equation}
\noindent
where $P_{m}(a)$ is a polynomial in $a$.  The integral above will be 
denoted by $N_{0,4}(a;m)$. This was the introduction of VHM to 
the wonderful world of definite integrals. George's original 
expression for the coefficients $d_{l}(m)$ in 
\begin{equation}
P_{m}(a) = \sum_{l=0}^{m} d_{l}(m) a^{l} 
\end{equation}
\noindent
was 
\begin{equation}
d_{l}(m) = \sum_{j=0}^{l} \sum_{s=0}^{m-l} \sum_{k=s+l}^{m} 
\frac{(-1)^{k-l-s}}{8^{k}} 
\binom{2k}{k} \binom{2m+1}{2s+2j} \binom{m-s-j}{m-k} \binom{s+j}{j} 
\binom{k-s-j}{l-j}. \nonumber 
\end{equation}
\noindent
Clearly not very useful. Intrigued by the fact that symbolic languages, such as
Mathematica or Maple, were unable 
to evaluate this integral, and being unaware of any useful technique available 
for the computation of integrals, the second author joined George Boros in 
the study of this polynomial. 

Computing some of these coefficients we observed that 
{\em they are all positive}. For example, if $m=5$, we obtain 
\begin{equation}
\{ d_{l}(5): \, 0 \leq l \leq 5 \}  = 
\left\{ \frac{4389}{256}, \, \frac{8589}{128}, \, \frac{7161}{64}, \, 
\frac{777}{8}, \, \frac{693}{16}, \, \frac{63}{8} \right\}.
\end{equation}
\noindent 
In this paper, we describe several 
proofs of the identity (\ref{int-01}) where the 
coefficients $d_{l}(m)$ are given by 
\begin{equation}
d_{l}(m) = 2^{-2m} \sum_{k=l}^{m} 2^{k} \binom{2m-2k}{m-k} \binom{m+k}{m} 
\binom{k}{l}.
\label{positive}
\end{equation}

\section{A survey of previous proofs}
\setcounter{equation}{0}

\subsection{The elementary proof} George's original idea was simple, but 
had profound consequences. The change of variables $x = \tan \theta$ 
yields 
\begin{equation}
N_{0,4}(a;m) = \int_{0}^{\pi/2} \left( \frac{\cos^{4} \theta}
{\sin^{4}\theta + 2a \sin^{2}\theta \cos^{2}\theta + \cos^{4}\theta} 
\right)^{m+1} \times 
\frac{d \theta}{\cos^{2} \theta}. \nonumber
\end{equation}
\noindent
Now, we observe that the denominator of the trigonometric function in the 
integrand, is a polynomial in $u = 2 \theta$. In detail, 
\begin{equation}
\sin^{4}\theta + 2a \sin^{2}\theta \cos^{2}\theta + \cos^{4}\theta 
 = (1+a) + (1-a) \cos^{2} u.
\nonumber
\end{equation}
\noindent 
In terms of the double-angle $u = 2 \theta$, the original integral becomes
\begin{equation}
N_{0,4}(a;m) = 2^{-(m+1)} \int_{0}^{\pi} 
\left( \frac{(1+ \cos u)^{2}}{(1+a) + (1-a) \cos^{2}u } \right)^{m+1} 
\times \frac{du}{1+ \cos u}. \nonumber
\end{equation}
\noindent
Expanding the binomial $(1+ \cos u)^{2m+1}$, the reader will check that 
\begin{equation}
\int_{0}^{\pi} \left[ (1 + a) + (1-a) \cos^{2}u \right]^{-(m+1)} 
\, \cos^{j} u \, du = 0, 
\label{vanishing}
\end{equation}
\noindent
for $j$ odd. The vanishing of half of the terms in the binomial expansion 
turned out to be a crucial property. The remaining integrals, those with 
$j$ even, can be simplified by using the double-angle trick one more time. 
The result is 
\begin{equation}
N_{0,4}(a;m)  =   \sum_{j=0}^{m} 2^{-j} \binom{2m+1}{2j} 
\int_{0}^{\pi} \left[ (3+a) + (1-a) \cos v \right]^{-(m+1)} ( 1 + \cos v)^{j} 
\, dv, \nonumber 
\end{equation}
\noindent
where $v = 2u$ and we have used the symmetry of cosine about $v = \pi$ to 
reduce the integrals form $[0, 2 \pi]$ to $[0, \pi]$. The  familiar 
change of variables $z = \tan(v/2)$ produces the  complicated formula for the 
coefficients $d_{l}(m)$ given in the Introduction. 

Excited with this evaluation and considering that we could not find it in
\cite{gr}, we wrote this as a paper and submitted it for publication. The 
reviewer was unforgiving: {\em Anybody who knows anything about hypergeometric
functions would be able to do this}.  We realized that it was time to 
learn about hypergeometric functions.  Eventually, we convinced an undergraduate
at Tulane University to check the details of the evaluation discussed 
above. They formed part of her final Honor Project and the results 
appeared  in \cite{sarah1}.

\subsection{The hypergeometric proof} 
The integral $N_{0,4}(a;m)$ can not be found in the compendium \cite{gr} in 
its original form. Naturally, it is possible that it appears 
in some disguised form, after a change of variables. In this case, the simplest 
change is $t = x^{2}$, that produces
\begin{equation}
N_{0,4}(a;m) = \frac{1}{2} \ift \frac{t^{-1/2} \, dt}{(t^{2} + 2at+1)^{m+1}}.
\label{nzero4}
\end{equation}
\noindent
Integrals of this type are listed in section $3.252$ of \cite{gr}. Indeed, 
$3.252.11$ states that 
\begin{equation}
\ift \frac{t^{\nu-1} \, dt}{(t^{2}+2at+1)^{\mu+1/2}}  = 
2^{-\mu} (a^{2}-1)^{\mu/2} \Gamma(1-\mu) B(\nu-2 \mu+1, -\nu) 
P_{\nu-\mu}^{\mu}(a). \nonumber 
\end{equation}
\noindent
The special functions involved here are the {\em gamma function} defined by
\begin{equation}
\Gamma(s) = \ift t^{s-1}e^{-t} \, dt,
\end{equation}
\noindent
that the reader must have seen as a generalization of integer 
factorials (namely that
$\Gamma(n) = (n-1)!$), and the {\em beta function} defined by
\begin{equation}
B(x,y) = \int_{0}^{1} t^{x-1}(1-t)^{y-1} \, dt,
\end{equation}
\noindent
that appears in beginning courses in Statistics, and finally the 
{\em associated Legendre function} $P_{\nu}^{\mu}(z)$ which was not very 
familiar to us. The table \cite{gr} is a good first source for special 
functions; in particular, sections $8.7$ and $8.8$ are 
dedicated to $P_{\nu}^{\mu}(z)$, 
where one learns that these functions 
are solutions of the differential equation
\begin{equation}
(1-z^{2}) \frac{d^{2}y}{dz^{2}} - 2z \frac{dy}{dz} + 
\left[ \nu (\nu+1) - \frac{\mu^{2}}{1-z^{2}} \right] y = 0. 
\end{equation}
\noindent
At the beginning of Section 8 in \cite{gr}, one finds the formula $8.702$:
\begin{equation}
P_{\nu}^{\mu}(z) = 
\frac{1}{\Gamma(1-\mu)} \left( \frac{z+1}{z-1} \right)^{\mu/2} 
{_{2}F_{1}} \left( -\nu, \nu+1; 1 - \mu; \, \frac{1-z}{2} \right). 
\end{equation}

The hypergeometric function $_{2}F_{1}(a,b;c;z)$ encountered above, 
is defined as the series
\begin{equation}
_{2}F_{1}(a,b;c;z) := \sum_{k=0}^{\infty} \frac{(a)_{k} \, (b)_{k}}
{(c)_{k} \, k!} z^{k},
\end{equation}
\noindent
where 
\begin{equation}
(a)_{k} = a(a+1)(a+2) \cdots (a+k-1) 
\end{equation}
\noindent
is the Pochhammer symbol. This function includes most of the elementary 
functions. For example, 
\begin{equation}
(1+z)^{n} = {_{2}F_{1}}(-n,b;b;-z), \quad \text{ for any }b,
\end{equation}
\noindent 
and 
\begin{equation}
\ln z = (z-1) {_{2}F_{1}}(1,1; \, 2, \, 1-z) 
\end{equation}
\noindent
just to name a couple. Section $9.121$ of \cite{gr} contain many more. 

Using these identities in (\ref{int-01}) produces
\begin{equation}
N_{0,4}(a;m) = 2^{m-1/2} (a+1)^{-(m+1/2)} B \left(2m+\tfrac{3}{2}, \tfrac{1}{2} 
\right) {_{2}F_{1}} \left[-m,m+1;m+ \tfrac{3}{2}; \, \tfrac{1-a}{2} \right],
\nonumber
\end{equation}
\noindent
and the formula (\ref{positive}) follows from here. The details appear in 
\cite{bomohyper}. The reviewer of \cite{sarah1} was correct: not much of 
the hypergeometric world is required to prove this. 

\subsection{ A detour into the world of Ramanujan} Not knowing about 
hypergeometric functions and before we search for this type of proof, we
observed that 
\begin{equation}
\ift \frac{dx}{bx^{4} + 2ax^{2}+1} = \frac{\pi}{2 \sqrt{2}}  
\frac{1}{\sqrt{a + \sqrt{b}}}. 
\label{int-1}
\end{equation}
\noindent
From here one sees that if we define 
\begin{equation}
g(c) := \ift \frac{dx}{x^{4} + 2ax^{2} + 1 + c} 
\end{equation}
\noindent
and 
\begin{equation}
h(c) = \sqrt{a + \sqrt{1+c}},
\end{equation}
\noindent
then (\ref{int-1}) shows that $g(c) = \pi \sqrt{2}h'(c)$. In particular,
\begin{equation}
h'(0) = \frac{1}{\pi \sqrt{2}} N_{0,4}(a;0). 
\label{der-0}
\end{equation}
\noindent
Now it is clear how to proceed: further differentiation yields 
\begin{equation}
\sqrt{a+ \sqrt{1+c}} = \sqrt{a+1} + \frac{1}{\pi \sqrt{2}} 
\sum_{k=1}^{\infty} \frac{(-1)^{k-1}}{k} N_{0,4}(a;k-1) c^{k}. 
\end{equation}
\noindent
Using Ramanujan Master's Theorem \cite{berndtI}, we were able to establish
the identity
\begin{equation}
\ift \frac{x^{m-1} \, dx}{(a + \sqrt{1+x})^{2m+1/2}} = 
\frac{1}{\pi} 2^{6m+3/2} \left[ m \binom{4m}{2m} \binom{2m}{m} \right]^{-1} 
\times N_{0,4}(a;m).  
\end{equation}
\noindent
The integral on the left can be evaluated in elementary form to produce 
(\ref{positive}). The details appear in \cite{bomoram}. It was never clear 
to the second  author, how George thought of using Ramanujan's Master 
Theorem. {\em He simply knew how to integrate}.

\subsection{Another elementary proof} The elementary proof described below
was shown to us by an (anonymous) referee. Write (\ref{der-0}) as 
\begin{equation}
\ift \frac{dx}{x^{4} + 2(b-1)x^{2} + 1+c} = \pi \sqrt{2} f'(c),
\end{equation}
\noindent
with 
\begin{equation}
f(c) = \sqrt{b + ( \sqrt{1+c} -1)}. 
\end{equation}
\noindent
Differentiating $m$ times with respect to $c$ yields 
\begin{equation}
m! (-1)^{m} \ift \frac{dx}{(x^{4}+2(b-1)x^{2} + 1 + c)^{m+1}}
= \pi \sqrt{2} f^{(m+1)}(c). 
\end{equation}
Hence
\begin{equation}
N_{0,4}(b-1;m) = \frac{\pi \sqrt{2}}{m!} (-1)^{m} f^{(m+1)}(0).
\end{equation}
\noindent
But 
\begin{equation}
f(c) = \sqrt{b} \left( 1 + \frac{\sqrt{1+c}-1}{b} \right)^{1/2} 
= \sqrt{b} \sum_{n=0}^{\infty} \binom{\tfrac{1}{2}}{n} b^{-n} ( \sqrt{1+c}-1
)^{n}. \nonumber
\end{equation}
\noindent
Now use the expansion
\begin{equation}
( \sqrt{1+c}-1)^{n} = \left( \frac{c}{2} \right)^{n} 
\left\{ 1 + n \sum_{i=1}^{\infty} \frac{(-1)^{i} (n+2i-1)!}{i! (n+i)!} 
\left( \frac{c}{4} \right)^{i} \right\}, \nonumber
\end{equation}
\noindent
given in \cite{gr}, Formula $1.114$, and the identity 
\begin{equation}
\binom{\tfrac{1}{2}}{n} = 
\frac{(-1)^{n-1}}{n2^{2n-1}} \binom{2n-2}{n-1}
\end{equation}
\noindent
and taking the coefficient of $c^{m+1}$ to obtain 
\begin{equation}
N_{0,4}(a;m) = \frac{\pi}{(8a+8)^{m+1/2}} \sum_{i=0}^{m} 
\binom{2m-2i}{k-i} \binom{m+i}{m} 2^{i} (a+1)^{i},
\end{equation}
\noindent
as desired. 

\section{A change of variables and the new proof}
\setcounter{equation}{0}

It soon became clear that the change of variables 
\begin{equation}
y = R_{2}(x) := \frac{x^{2}-1}{2x} 
\label{transf-r2}
\end{equation}
\noindent
was the key to extend the elementary method to integrals of higher order. The
inverse has two branches 
\begin{equation}
x = y \pm \sqrt{y^{2}+1},
\end{equation}
\noindent
where the plus sign is valid for $x \in [0, \, +\infty)$ and the other one 
on $(-\infty, 0]$. The rational function $R_{2}$ arises from the 
identity
\begin{equation}
\cot 2 \theta = R_{2}(\cot \theta). 
\end{equation}
\noindent
This change of variables gives the proof of the next theorem and it 
leads to interesting equality among integrals.

\begin{Thm}
\label{thm-inv}
Let $f$ be a rational function and assume that the integral of $f$ over 
$\mathbb{R}$ is finite. Then 
\begin{eqnarray}
\ifft f(x) \, dx & =  & 
\ifft \left[ f(y + \sqrt{y^{2}+1}) + f(y - \sqrt{y^{2}+1}) 
\right] \, dy +  \label{invar} \\
& + & \ifft \left[ f(y + \sqrt{y^{2}+1}) - f(y - \sqrt{y^{2}+1}) 
\right] \, \frac{y \, dy}{\sqrt{y^{2}+1}}.
\nonumber
\end{eqnarray}
\noindent
Moreover, if $f$ is an 
{\em even} rational function,  the identity (\ref{invar}) remains valid by
replacing the interval of integration from $\mathbb{R}$ to ${\mathbb{R}}^{+}$.
\end{Thm}

\medskip

The reader will check that the integrand on the right hand side of 
(\ref{invar}) is also a rational function. This new function is an 
even rational function if $f$ is such. 

The first
application of Theorem \ref{thm-inv} was given in \cite{boros2} where we 
proved that the integral 
\begin{equation}
U_{6}(a,b;c,d,e) := \ift \frac{cx^{4}+dx^{2}+e}{x^{6}+ax^{4}+bx^{2}+1} \, dx
\end{equation}
\noindent
is invariant under the change of parameters 
\begin{eqnarray}
a & \to & \frac{ab + 5a + 5b + 9}
{( a + b + 2)^{4/3} } \label{scheme6a}  \\
b & \to & \frac{a + b+6}{( a + b + 2)^{2/3} } \nn \\
c & \to & \frac{c + d + e}{(a+b+2)^{2/3} } \nn \\ 
d & \to & \frac{c(b + 2) + 2d + e(a + 3)}{ a + b + 2} \nn \\
e & \to & \frac{c + e}{( a+ b + 2)^{1/3} }. \nn 
\end{eqnarray}
\noindent
Iteration of (\ref{scheme6a}) produces a sequence of parameters 
$(a_{n},b_{n};c_{n},d_{n},e_{n})$ with 
\begin{equation}
U_{6}(a_{n},b_{n};c_{n},d_{n},e_{n}) = 
U_{6}(a,b;c,d,e). \label{invar-int}
\end{equation}
\noindent
Moreover, as $n \to \infty$, we have that $a_{n}, \, b_{n} \to 3$ and there 
exists a number $L$ such that $(c_{n}, d_{n}, e_{n}) \to L(1,2,1)$. Passing to 
the limit in (\ref{invar-int}) we obtain
\begin{equation}
U_{6}(a,b;c,d,e) = \frac{\pi}{2}L. 
\label{limit-rat}
\end{equation}
\noindent
This is precisely the rational function analog of the transformation of 
parameters 
\begin{eqnarray}
a & \mapsto & \frac{a+b}{2} \label{elliptic} \\
b & \mapsto & \sqrt{ab}. \nonumber 
\end{eqnarray}
\noindent
that leaves the elliptic integral 
\begin{equation}
G(a,b) = \int_{0}^{\pi/2} \frac{d \varphi}{\sqrt{a^{2} \cos^{2} \varphi + 
b^{2} \sin^{2} \varphi}}
\end{equation}
\noindent
invariant. In this 
situation, the sequence $(a_{n},b_{n})$ obtained by iterating 
(\ref{elliptic}) converges to a common limit: the {\em arithmetic-geometric
mean} of $a$ and $b$, denoted by $\text{AGM}(a,b)$. In the limit, we have 
\begin{equation}
G(a,b) = \frac{\pi}{2 \text{AGM}(a,b)}, 
\end{equation}
\noindent
precisely as in (\ref{limit-rat}).  Information about the $\text{AGM}$
can be found in \cite{borwein1} and \cite{cox85}.

This was the beginning of a series of results on changes of parameters that 
leave integrals of rational functions invariant. We have called 
(\ref{scheme6a}) a {\em rational Landen transformation}. These were 
extended to all 
even rational functions in \cite{boros1} and the convergence of the iterations 
of these transformation were described in \cite{hubbard1}. The extension 
to {\em arbitrary} rational functions
on the whole line is described in  
\cite{manna-moll2} and \cite{manna-moll1}. It is an open question to 
produce rational Landen transformations for a general rational function on 
the interval $[0, \infty)$. The vanishing of the integrals such as 
(\ref{vanishing}) does not happen if the integrand is not an even rational 
function. 

From the 
point of view of numerical calculations, the rational Landen transformations
give a procedure to evaluate definite integrals in iterative form. For 
example, in the simplest case of  the elementary integral
\begin{equation}
I = \ifft \frac{dx}{ax^{2}+bx+c},
\end{equation}
\noindent
with $b^{2}-4ac < 0$, we have shown in \cite{manna-moll1} that $I$ is 
invariant under the change of parameters
\begin{eqnarray}
a_{n+1} & = & a_{n} \left[ \frac{(a_{n}+3c_{n})^{2}-3b_{n}^{2}}
{(3a_{n}+c_{n})(a_{n}+3c_{n})-b_{n}^{2}} \right], \label{cubic} \\
b_{n+1} & = & b_{n} \left[ \frac{3(a_{n}-c_{n})^{2}-b_{n}^{2}}
{(3a_{n}+c_{n})(a_{n}+3c_{n})-b_{n}^{2}} \right], \nonumber  \\
c_{n+1} & = & c_{n} \left[ \frac{(3a_{n}+c_{n})^{2}-3b_{n}^{2}}
{(3a_{n}+c_{n})(a_{n}+3c_{n})-b_{n}^{2}} \right]. \nonumber 
\end{eqnarray}
\noindent
Moreover, the error defined by 
\begin{equation}
e_{n} = (a_{n} - \tfrac{1}{2} \sqrt{4ac-b^{2}}, b_{n}, c_{n} - 
\tfrac{1}{2} \sqrt{4ac-b^{2}} ), 
\end{equation}
\noindent
decays to $0$ at a cubic rate, that is, $\vert| e_{n+1} \vert| 
\leq C \vert| e_{n} \vert|^{3}$. 

\medskip

The transformation (\ref{transf-r2}) and its higher order analogues, given by
the rational function $R_{m}$ defined by
\begin{equation}
\cot(m \theta) = R_{m}(\cot \theta),
\end{equation}
\noindent
have produced interesting results for definite integrals of higher degree. 
See \cite{manna-moll2} for details. 

The second goal of
this note is to use the Landen transformation (\ref{transf-r2}) 
to produce a new proof of the identity (\ref{int-01}) and thus obtain the 
more explicit natural form of the coefficients $d_{l}(m)$ given in 
(\ref{positive}). 

\begin{Thm}
For $m \in \mathbb{N}$,  let 
\begin{equation}
Q(x) = \frac{1}{(x^{4}+2ax^{2}+1)^{m+1}}. 
\end{equation}
\noindent
Define
\begin{eqnarray}
Q_{1}(y) & := & \left[ Q(y+\sqrt{y^{2}+1}) + Q(y-\sqrt{y^{2}+1}) 
\right] +  \nonumber \\
& + & \frac{y}{\sqrt{y^{2}+1}} 
\left[ Q(y+\sqrt{y^{2}+1}) - Q(y-\sqrt{y^{2}+1}) 
\right].   \nonumber
\end{eqnarray}
\noindent
Then 
\begin{equation}
Q_{1}(y) =  \frac{T_{m}(2y)}{2^{m}(1+a+2y^{2})^{m+1}},
\end{equation}
\noindent
where 
\begin{equation}
T_{m}(y) = \sum_{k=0}^{m} \binom{m+k}{m-k} y^{2k}. 
\label{bin-sum}
\end{equation}
\end{Thm}
\begin{proof}
Introduce the variable  $\phi = y + \sqrt{y^{2}+1}$. Then $y - \sqrt{y^{2}+1} = 
- \phi^{-1}$ and  $y = \phi - \phi^{-1}$.  Moreover, 
\begin{eqnarray}
Q_{1}(y) & = & \left[ Q(\phi) + Q(\phi^{-1}) \right] + 
\frac{\phi^{2}-1}{\phi^{2}+1} 
\left( Q(\phi) - Q(\phi^{-1}) \right) \nonumber \\
 & = & \frac{2}{\phi^{2}+1} \left[ \phi^{2}Q(\phi) + Q(\phi^{-1}) 
\right] \nonumber \\
& := & S_{m}(\phi). \nonumber 
\end{eqnarray}
\noindent
The result of the theorem is equivalent to
\begin{equation}
2^{m} \left(1 + a + \tfrac{1}{2}(\phi - \phi^{-1})^{2} ) \right)^{m+1} \, 
S_{m}(\phi) = T_{m}(\phi - \phi^{-1}). 
\label{newform}
\end{equation}
\noindent
A direct simplification of the left hand side of (\ref{newform}) shows that 
this identity is equivalent to proving
\begin{equation}
\frac{\phi^{2m+1}+ \phi^{-(2m+1)}}{\phi+\phi^{-1}} = T_{m}(\phi - \phi^{-1}). 
\label{newform2}
\noindent
\end{equation}

\noindent
{\bf First proof}. One 
simply checks that both sides of (\ref{newform2}) satisfy the 
second order recurrence 
\begin{equation}
c_{m+2}-  ( \phi^{2} + \phi^{-2}) c_{m+1} + c_{m} = 0, 
\label{rec-22}
\end{equation}
\noindent
and the values for $m=0$ and $m=1$ match. This is straight-forward for the 
expression on the left hand side, while the WZ-method 
settles the right hand side. 
\end{proof}

\medskip

\no
{\bf Note}. The WZ-method is an algorithm  in Computational Algebra that, among 
other things, will produce for a hypergeometric/holonomic 
sum, such as (\ref{bin-sum}), a 
recurrence like (\ref{rec-22}). The reader will find in \cite{aequalsb} 
and \cite{nemes} information about this algorithm.  \\

\noindent
{\bf Second proof}. In \cite{graham1}, one finds the generating function
\begin{equation}
B_{t}(z) = \sum_{k \geq 0} (tk)_{k-1} \frac{z^{k}}{k!}, 
\end{equation}
\noindent
where $(a)_{k} = a(a+1) \ldots (a+k-1)$ is the Pochhammer symbol. The 
special values
\begin{equation}
B_{-1}(z) = \frac{1+ \sqrt{1+4z}}{2}, \quad \text{ and } 
B_{2}(z) = \frac{1- \sqrt{1-4z}}{2z}, 
\end{equation}
\noindent
are combined to produce the identity
\begin{equation}
\frac{1}{\sqrt{1+4z}} \left( B_{-1}(z)^{n+1} - (-z)^{n+1} B_{2}(-z)^{n+1} 
\right) = \sum_{k=0}^{n} \binom{n-k}{k} z^{k}. 
\end{equation}
\noindent
Replace $n$ by $2m$ and $z$ by $(4y^{2})^{-1}$ to produce 
\begin{equation}
\frac{1}{2 \sqrt{1+y^{2}} (2y)^{2m}} ( \phi^{2m+1} + \phi^{-(2m+1)} ) = 
\sum_{k=0}^{m} \binom{2m-k}{k} z^{k}.
\nonumber
\end{equation}
\noindent
We also have 
\begin{eqnarray}
T_{m}(y)  & = & \sum_{k=0}^{m} \binom{m+k}{m-k} y^{2k} \nonumber \\
 & = & y^{2m} \sum_{k=0}^{m} \binom{2m-k}{k} y^{-2k}. \nonumber 
\end{eqnarray}
\noindent
Thus,
\begin{eqnarray}
T_{m}(\phi - \phi^{-1}) & = & T_{m}(2y) \nonumber \\
 & = & (2y)^{2m} \sum_{k=0}^{m} \binom{2m-k}{k} z^{k}. \nonumber
\end{eqnarray}
\noindent
We conclude that 
\begin{equation}
T_{m}(\phi - \phi^{-1}) = \frac{1}{2 \sqrt{1+y^{2}}} (\phi^{2m+1} + 
\phi^{-(2m+1)} ), 
\end{equation}
\noindent
and the result follows from $\phi + \phi^{-1}  = 2 \sqrt{y^{2}+1}$.

\bigskip

We now prove (\ref{int-01}). 
The identity in Theorem 3.1 shows that 
\begin{equation}
\ift Q(x) \, dx  =  \ift Q_{1}(y) \, dy, 
\label{equalint}
\end{equation}
\noindent
and this last integral can be evaluated in elementary terms. Indeed, 
\begin{eqnarray}
\ift Q_{1}(y) \, dy & = & \ift \frac{T_{m}(2y) \, dy}{2^{m} 
(1+ 2y^{2})^{m+1}} \nonumber \\
& = & \frac{1}{2^{m}} \sum_{k=0}^{m} \binom{m+k}{m-k} 
\ift \frac{(2y)^{2k} \, dy}{(1 + a + 2y^{2})^{m+1}}. \nonumber
\end{eqnarray}
\noindent
The change of variables $y = \frac{\sqrt{1+a}}{\sqrt{2}} \, t$ gives
\begin{equation}
\ift Q_{1}(y) \, dy = \frac{1}{[2(1+a)]^{m+1/2}} 
\sum_{k=0}^{m} \binom{m+k}{m-k} 2^{k} (1+a)^{k} 
\ift \frac{t^{2k} \, dt}{(1+t^{2})^{m+1}}, \nonumber 
\end{equation}
\noindent
and the elementary identity
\begin{equation}
\ift \frac{t^{2k} \, dt}{(1+t^{2})^{m+1}} = \frac{\pi}{2^{2m+1}} 
\binom{2k}{k} \binom{2m-2k}{m-k} \binom{m}{k}^{-1}
\nonumber
\end{equation}
gives
\begin{equation}
\ift Q_{1}(y) \, dy = \frac{\pi}{2^{2m+1}} \frac{1}{[2(1+a)]^{m+1/2}} 
\sum_{k=0}^{m} \binom{m+k}{m-k} 2^{k} 
\binom{2k}{k} \binom{2m-2k}{m-k} \binom{m}{k}^{-1}
(1+a)^{k}.
\nonumber
\end{equation}
\noindent
This can be simplified further using 
\begin{equation}
\binom{m+k}{m-k} \binom{2k}{k} = \binom{m+k}{m} \binom{m}{k},
\end{equation}
\noindent 
and the equality (\ref{equalint}) to produce 
\begin{equation}
\ift Q(y) \, dy = \frac{\pi}{2^{2m+1}} \frac{1}{[2(1+a)]^{m+1/2}} 
\sum_{k=0}^{m} 2^{k} \binom{m+k}{m}  \binom{2m-2k}{m-k} (1+a)^{k}. 
\end{equation}
This completes the proof of (\ref{int-01}),  with 
the coefficients $d_{l}(m)$ given in 
(\ref{positive}).

\bigskip

\no
{\bf Acknowledgments}. The work of the second author was partially funded by
$\text{NSF-DMS } 0409968$.

\end{document}